\documentclass{amsart}
\usepackage{a4wide, graphicx, amsmath, amsfonts, amsthm, amssymb, latexsym, tikz, color, enumitem,  verbatim, longtable}   

\usepackage[colorlinks]{hyperref}
\hypersetup{linkcolor=blue, urlcolor=blue, citecolor=red}

\let\oldmarginpar\marginpar
\renewcommand\marginpar[1]{\-\oldmarginpar[\raggedleft\footnotesize #1]%
{\raggedright\footnotesize #1}}

\newtheorem{thm}{Theorem}[section]
\newtheorem{lem}[thm]{Lemma}
\newtheorem{prop}[thm]{Proposition}

\newtheorem{defn}[thm]{Definition}

\newcommand{\set}[2]{\{#1:#2\}}
\newcommand{\genset}[1]{\langle#1\rangle}

\newcommand{\sym}{\operatorname{Sym}(\Omega)}

\newcommand{\N}{\mathbb{N}}

\newcommand{\all}{S_{1,2,3,4,5}}
\renewcommand{\and}{\text{ and }}

\newcommand{\T}{\Omega^\Omega}

\author{J. Jonu\v{s}as and J. D. Mitchell}
\title[A finite interval in the subsemigroup lattice of $\Omega^\Omega$]{A finite interval in the subsemigroup lattice of the full transformation monoid}
\date{\today}

\begin{document}
 
\begin{abstract}
In this paper we describe a portion of the subsemigroup lattice of the \emph{full transformation semigroup} $\Omega^\Omega$,
 which consists of all mappings on the 
countable infinite set $\Omega$. 
Gavrilov showed that there are five maximal subsemigroups of $\Omega^\Omega$ containing the symmetric group $\sym$. 
The portion of the subsemigroup lattice of  $\Omega^\Omega$ which we describe is that between the intersection of these five maximal subsemigroups and $\Omega^\Omega$. 
We prove that there are only 38 subsemigroups in this interval, in contrast to the $2^{2^{\aleph_0}}$ subsemigroups between $\sym$ and $\Omega^\Omega$.
\end{abstract}

\maketitle

\section{Introduction}

Let $\Omega$ denote an arbitrary  infinite set, let $\Omega^\Omega$ denote the 
semigroup of mappings from $\Omega$ to itself, and let $\sym$ denote the symmetric group on $\Omega$.
The subsemigroups of $\Omega^\Omega$ form an algebraic lattice with $2^{|\Omega|}$ compact (finitely generated) 
elements under inclusion. 
Pinsker and Shelah \cite{Pinsker2011aa} proved that every algebraic lattice with at most $2^{|\Omega|}$ compact 
elements 
can be embedded into the subsemigroup lattice of $\Omega^\Omega$. There are $2^{2^{|\Omega|}}$ distinct 
subsemigroups of $\Omega^\Omega$, and even this many maximal subsemigroups (see, for example, \cite{Goldstern2002aa} 
or \cite{Rosenberg1976aa}). In contrast, 
in 1965 Gavrilov \cite{Gavrilov1965aa} showed that there are five maximal subsemigroups of $\Omega^\Omega$ containing 
$\sym$
when $\Omega$ is countable. In  2005, Pinsker \cite{Pinsker2005aa} extended Gavrilov's result to sets of arbitrary infinite cardinality 
showing that there are $2|\alpha|+1$ such maximal subsemigroups when $|\Omega|=\aleph_{\alpha}$. 
Clone lattices on infinite sets have also been extensively studied, see \cite{Beiglbock2009aa} and \cite{Goldstern2008aa}. 

If $S$ is a subsemigroup of a semigroup $T$, then  we denote this by $S\leq T$. 
If $S$ and $T$ are subsemigroups of $\T$, then the \textit{interval from $T$ to $S$}, denoted $(T,S)$, is defined to be
the set of proper subsemigroups of $S$ properly containing $T$, i.e.
$$(T,S) = \set{U\leq \Omega^\Omega}{T \lneq U \lneq S }.$$

If $S$ is any maximal subsemigroup of $\T$ containing $\sym$, as described by Gavrilov \cite{Gavrilov1965aa} and Pinsker 
\cite{Pinsker2005aa}, then the 
interval $(\sym, S)$ has cardinality $2^{2^{\kappa}}$ where  $|\Omega|=\aleph_{\alpha}$ and 
$\kappa=\max\{\alpha, \aleph_0\}$;  for further details see Pinsker \cite{Pinsker2005ab}.  
The  maximal subsemigroups of the maximal subsemigroups described by Gavrilov \cite{Gavrilov1965aa} are classified 
in \cite{Dougall2012aa}; perhaps 
surprisingly there are only countably many such semigroups. In further contrast to Pinsker and Shelah's result 
\cite{Pinsker2011aa}, 
we prove that there are only 36 subsemigroups in the interval from 
the intersection of the maximal subsemigroups described by Gavrilov \cite{Gavrilov1965aa} to $\Omega^\Omega$ when 
$\Omega$ is countably infinite; we completely describe the subsemigroup lattice in this interval.

Throughout this paper, we write functions to the right of their argument and compose from left to right. 
If $\alpha\in \Omega$, $f\in\Omega^\Omega$ and $\Sigma\subseteq \Omega$, then  
$\alpha f^{-1}=\set{\beta\in \Omega}{\beta f=\alpha}$, $\Sigma f=\set{\alpha f}{\alpha\in \Sigma}$, and $f|_{\Sigma}$ 
denotes the \emph{restriction} of $f$ to $\Sigma$. 
If $f\in\Omega^\Omega$ and $\Sigma\subseteq \Omega$ such that  $f|_{\Sigma}$ is injective and $\Sigma f=\Omega f$, 
then we will refer to $\Sigma$ as a \emph{\hypertarget{transversal}{transversal}} of $f$. 

We require the following parameters of a function $f\in \Omega^\Omega$ to  describe the maximal subsemigroups 
containing  $\sym$ and the semigroups introduced here:
\begin{eqnarray*}
d(f)& =& |\Omega\setminus\Omega f|\\
c(f)& =&|\Omega\setminus \Sigma|\text{, where }\Sigma\text{ is any transversal of }f\\
k(f)& =& |\set{\alpha\in\Omega}{|\alpha f^{-1}|=\infty}|.
\end{eqnarray*}
 The parameters $d(f),c(f)$, and $k(f)$ were termed the \emph{defect}, \emph{collapse}, and \emph{infinite 
 contraction index}, respectively, of $f$ in~\cite{Howie1998aa}. Recall that the \emph{kernel} of $f\in \Omega^\Omega$ is the equivalence relation:
 $$\ker(f)=\set{(\alpha, \beta)\in \Omega\times \Omega}{\alpha f=\beta f}.$$
If $S$ is a semigroup and $T$ is a subset of $S$, then by $\langle T \rangle$ we 
denote the subsemigroup of $S$ generated by $T$. 

Throughout the remainder of the paper we will denote by $\Omega$ an arbitrary countably infinite set. 
Using the notation of \cite{east}, the maximal semigroups of $\T$ containing $\sym$ are:
\begin{align*}
  S_1 =& \{ f \in \T : c(f) = 0 \text{ or } d(f) > 0 \} \\
  S_2 =& \{ f \in \T : c(f) > 0 \text{ or } d(f) = 0 \} \\
  S_3 =& \{ f \in \T : c(f) < \infty \text{ or } d(f) = \infty \} \\
  S_4 =& \{ f \in \T : c(f) = \infty \text{ or } d(f) < \infty \} \\
  S_5 =& \{ f \in \T : k(f) < \infty \}. 
\end{align*}
For the sake of brevity, if $I\subseteq \{1,2,3,4,5\}$, then we will denote the intersection $(\bigcap_{i\in I} S_i)$ 
by $S_I$, so that $S_{1,2}$ denotes $S_1\cap S_2$ and so on. \vspace{\baselineskip}

\noindent{\bf Main Theorem.} {\it
The interval $(\all, \Omega^\Omega)$ consists of the $30$ semigroups which are all possible intersections of any non-empty proper subset of 
$\{S_1, S_2, S_4, S_3, S_5\}$  and the semigroups: 
$$U,\ V,\ S_1\cap U,\ S_2\cap V,\ S_5\cap U,\  S_5\cap V,\ S_{1,5}\cap U, \text{ and  }S_{2,5}\cap V,$$
where
\begin{eqnarray*}
U&=&\{ f \in \T : d(f) = \infty \text{ or } (0<c(f)<\infty) \}\cup \sym\\
V&= &          \{f \in \T : c(f)=\infty \text{ or } (0<d(f)<\infty) \} \cup \sym.
\end{eqnarray*}}

We will describe the relevant semigroups given in the Main Theorem, in terms of the parameters $d$, $c$, and $k$, in Section 
\ref{description}.
We will show in Proposition \ref{U_are_a_semigroup} that $U$ and $V$ are semigroups. It is routine to verify that 
$U$ is contained in $S_3$ and 
$S_2 \cap S_3 \leq U$.
But there exists $f\in U$ such that $d(f)=\infty$ and $c(f)=0$ and there is no such function in 
$S_2\cap S_3$, and so
$S_2\cap S_3\not=U$. Similarly, $V$ is contained in $S_4$ and $S_2 \cap S_4 \lneq V$.

In Section \ref{sec:lemmas} a 
number of technical lemmas are proved which are extensively used in the proof of the Main 
Theorem in Section 
\ref{sec:proof}. 


\section{Descriptions of the semigroups in the Main Theorem}\label{description}

In this section we describe the semigroups in the Main Theorem according to the parameters 
$c$, $d$, and $k$.  We will not repeat these statements later in the paper, but have opted to collect the descriptions in one place, for
ease of reference. 

\begin{center}
\begin{longtable}{|c|l|}\hline
$S$ &Description of the elements $f\in S$\\
\hline\hline
 $S_1 $&$ c(f) = 0 \text{ or } d(f) > 0$\\\hline
 $S_2 $&$ c(f) > 0 \text{ or } d(f) = 0$\\\hline
 $S_3 $&$ c(f) < \infty \text{ or } d(f) = \infty$\\\hline
 $S_4 $&$ c(f) = \infty \text{ or } d(f) < \infty$\\\hline
 $S_5 $&$ k(f) < \infty $\\\hline\hline

$U $&$ d(f) = \infty \text{ or } (0<c(f)<\infty) \text{ or } f\in \sym$\\\hline
 $V $&$ c(f) = \infty \text{ or } (0<d(f)<\infty) \text{ or } f\in \sym$\\\hline\hline

 $S_{1,2} $&$ c(f)>0 \and d(f)>0 \text{ or } f\in \sym$\\\hline
 $S_{1,3} $&$ d(f) = \infty \text{ or } c(f)=0 \text{ or } (0<c(f),d(f)<\infty)$\\\hline
 $S_{1,4} $&$ (c(f)=d(f) = \infty) \text{ or } (0<d(f)<\infty) \text{ or }  f\in \sym$\\\hline
 $S_{1,5} $&$ (k(f) <\infty \and  d(f)>0) \text{ or } c(f) =0$\\\hline
 $S_{2,3} $&$ (c(f)=d(f) = \infty) \text{ or } (0<c(f)<\infty) \text{ or } f\in \sym$\\\hline
 $S_{2,4} $&$ c(f) = \infty \text{ or } d(f)=0 \text{ or } (0<c(f),d(f)<\infty)$\\\hline
 $S_{2,5} $&$ (k(f) <\infty \and  c(f)>0) \text{ or } (d(f) =0 \and k(f)<\infty)$\\\hline
 $S_{3,4} $&$ (c(f)=d(f)=\infty) \text{ or }(c(f), d(f)< \infty)  $\\\hline
 $S_{3,5} $&$ c(f)<\infty \text{ or } (d(f)=\infty \and k(f)<\infty) $\\\hline
 $S_{4,5} $&$ (c(f)=\infty \and k(f)<\infty) \text{ or } (d(f)<\infty \and k(f)<\infty)$\\\hline\hline

 $V \cap S_2 $&$ c(f) = \infty \text{ or } (0 < c(f),d(f) <\infty) \text{ or }  f\in \sym$\\\hline
 $U \cap S_1 $&$ d(f) = \infty \text{ or } (0 < c(f),d(f) <\infty) \text{ or }  f\in \sym$\\\hline
 $V \cap S_5 $&$ (k(f)<\infty \and c(f)= \infty) \text{ or } (0 < d(f)<\infty \and k(f)<\infty) \text{ or } f\in \sym$\\\hline
 $U \cap S_5 $&$ (k(f)<\infty \and d(f)= \infty) \text{ or } (0 < c(f)<\infty) \text{ or } f\in \sym$\\\hline\hline

 $S_{1,2,3} $&$ (c(f)>0 \and d(f) = \infty) \text{ or } (0 < c(f),d(f) <\infty) \text{ or } f\in \sym$\\\hline
 $S_{1,2,4} $&$ (c(f)=\infty \and d(f) > 0) \text{ or } (0 < c(f),d(f) <\infty) \text{ or } f\in \sym$\\\hline
 $S_{1,2,5} $&$ (c(f), d(f)>0) \and k(f)< \infty \text{ or } f\in \sym$\\\hline
 $S_{1,3,4} $&$ (c(f)=d(f)= \infty) \text{ or } (0<d(f)<\infty \and  c(f)<\infty) \text{ or } f\in \sym$\\\hline
 $S_{1,3,5} $&$ (d(f)= \infty \and k(f)<\infty) \text{ or } (c(f)<\infty  \and  d(f)>0) \text{ or } f\in \sym$\\\hline
 $S_{1,4,5} $&$ (c(f)=d(f)= \infty \and k(f)<\infty) \text{ or }  (0 < d(f)<\infty \and k(f)<\infty) \text{ or } f\in \sym$\\\hline
 $S_{2,3,4} $&$ (c(f)=d(f)= \infty) \text{ or } ( 0 < c(f) <\infty \and d(f)<\infty)\text{ or } f\in \sym$\\\hline
 $S_{2,3,5} $&$ (c(f)=d(f)= \infty \and k(f)<\infty) \text{ or }  (0 <c(f)<\infty)\text{ or } f\in \sym$\\\hline
 $S_{2,4,5} $&$ (c(f)=\infty \and k(f)<\infty) \text{ or } ( 0 <c(f) <\infty \and d(f)<\infty)\text{ or } f\in \sym$\\\hline
 $S_{3,4,5} $&$ (k(f)<\infty \and c(f)=d(f)=\infty )\text{ or } (c(f),d(f) < \infty ) $\\\hline\hline

 $V \cap S_{2,5}$&$ (k(f)<\infty \and c(f)=\infty) \text{ or } (0<c(f), d(f)<\infty)\text{ or } f\in \sym$\\\hline
 $U \cap S_{1,5}$&$ (k(f)<\infty \and d(f)=\infty) \text{ or } (0<c(f), d(f)<\infty)\text{ or } f\in \sym$\\\hline\hline

 $S_{1,2,3,4} $&$ (c(f)=d(f)= \infty) \text{ or } (0 < c(f), d(f)<\infty)  \text{ or } f\in \sym $\\\hline
 $S_{1,2,3,5} $&$ (k(f)<\infty , c(f)> 0 \and d(f) = \infty) \text{ or } (0< c(f), d(f)<\infty) \text{ or } f\in \sym$\\\hline
 $S_{1,2,4,5} $&$ k(f)< \infty \and  [(c(f)=d(f)=\infty) \text{ or }(c(f)>0  \and 0< d(f) < \infty)] \text{ or } f\in \sym$\\\hline
 $S_{1,3,4,5} $&$ (k(f)< \infty \and  c(f)=d(f)=\infty) \text{ or }(c(f)<\infty  \and 0< d(f) < \infty) \text{ or } f\in \sym$\\\hline 
 $S_{2,3,4,5} $&$ (k(f)< \infty \and  c(f)=d(f)=\infty) \text{ or }(d(f)<\infty  \and 0< c(f) < \infty) \text{ or } f\in \sym$\\\hline\hline

 $\all $&$ (c(f)=d(f)=\infty \and k(f)<\infty) \text{ or } (0<c(f),d(f)<\infty)\text{ or } f\in \sym$\\\hline
\end{longtable} 
\end{center}


\section{Technical lemmas} \label{sec:lemmas}

We require several technical results to prove the Main Theorem, which we present in this section.
We will make repeated use of the following lemma without reference. 

\begin{lem}[Lemma 5.4 in \cite{east}] \label{inequalities}
Let $u,v \in \T$. Then the following hold:
\begin{enumerate} [label = \rm{(\roman*})]
 \item $d(g)\leq d(fg)\leq d(f)+d(g)$;
 \item   $c(v) \leq c(vu) \leq c(v)+c(u)$;
 \item if $c(g)<\infty=d(f)$, then $d(fg)=\infty$;
 \item  if $d(f)<\infty=c(g)$, then $c(fg)=\infty$;
 \item $k(fg)\leq k(f)+k(g)$.
\end{enumerate}
\end{lem}

We will repeatedly use  the observation
that if $f,g\in \Omega^\Omega$ such that $\ker(f)=\ker(g)$ and $d(f)=d(g)$, then $f\in \genset{\sym, g}$. 
Note that if $m,n\in \N\cup \{\infty\}$, then there exists $f\in \Omega^\Omega$ such that $c(f)=m$ and $d(f)=n$. 
Also note that $f\in \Omega^\Omega$ is injective if and only if $c(f)=0$ and $f$ is surjective if and only if $d(f)=0$. 

\begin{lem} \label{tech}
Let $u,v, f\in \T$. Then the following hold:
\begin{enumerate} [label = \rm{(\roman*})]
\item \label{c finite}
if $0<c(u), c(f)<\infty$ and $d(u)=d(f)$, then  $f\in\genset{\all, u}$;
\item \label{c finite more}
if $c(u)<\infty$, $d(u)=\infty$, $c(v)>0$, $d(f)=d(v)=0$, and $0<c(f)<\infty$, then 
$f\in  \genset{\all, u,v}$;
\item \label{d infinite} if $u$ is injective, $c(f)<\infty$, and $d(u)=d(f)=\infty$, then $f\in \genset{\all, u}$;
\item \label{two fun} if $u$ is injective, $d(u) > 0$, $d(v)=d(f)=\infty$, and $c(v), c(f)<\infty$ ,
then $f \in \genset{\all, u, v}$.
\end{enumerate}
\end{lem}
\proof 
{\bf (i).} Let $g\in \T$ be such that $\ker(g) = \ker(f)$ and such that  $ \Omega g$ is a  transversal of $u$. Then $c(g) = c(f)$ 
and $d(g) = c(u)$ and  thus $g \in \all$. Since $d(gu)=d(u)=d(f)$ and $\ker(gu)=\ker(g)=\ker(f)$,  $f\in \genset{gu, \sym}\leq \genset{u, \all}$.
\vspace{\baselineskip}

\noindent{\bf (ii).} If $c(v) <\infty$, then $f \in \genset{\all, v}$ by Lemma \ref{tech}\ref{c finite}. Suppose $c(v)=\infty$ and let $t \in \all$ 
be such that $0<c(t),d(t)<\infty$. Then $0<c(tu)<\infty$ and $d(tu)=\infty$. Choose $g \in \T$ 
such that $\Omega tu$ is a transversal of $g$, $\Omega g$ is a transversal of $v$ and $k(g)<\infty$. Then $c(g)=d(g)=\infty$ and 
so $g\in \all$.
Also $0<c(tugv)<\infty$ and $d(tugv)=0$ and thus $f \in \genset{\all,tugv} \leq \genset{\all, u, v}$ by Lemma \ref{tech}\ref{c finite}.
\vspace{\baselineskip}

\noindent{\bf (iii).} If $c(f)=0$, then $\ker(f)=\ker(u)$ and so $f\in \genset{u, \all}$. 
  
  Suppose $0<c(f)<\infty$. If $g \in \T$ is such that $\ker(g) = \ker(f)$ and $0<d(g)<\infty$, then $g\in \all$.  
But $\ker(gu)=\ker(g)$ since $u$ is injective and so $\ker(gu)=\ker(f)$. Also $d(gu)\geq d(u)=\infty$, and so $d(gu)=d(u)=d(f)$. Thus  $f\in  \genset{gu, \sym}\leq \genset{u, \all}$. \vspace{\baselineskip}

\noindent{\bf (iv).} If $v$ is injective or $d(u)=\infty$, then $f\in \genset{\all, v}$ or $f\in \genset{\all, u}$, respectively, by part \ref{d infinite}. 
Suppose $0<c(v)$ and  $d(u) < \infty$. 
 Let $w \in \T$ be such that $c(w)=0$ and $\Omega w$ is a transversal of $v$. Then $c(wv)=0$ and $d(wv)>d(v)=\infty$, 
 and so  $f \in \genset{\all, wv}$ by part \ref{d infinite}. 
Let $t \in \T$ be any function such that $\Omega u$ is a transversal of $t$ and $d(t)=d(w)$. Then $c(t)=d(u)$ and $0<c(t)<\infty$. 
Since $d(w)=c(v)$, it follows that $0<d(t)<\infty$. Hence $t\in \all$.  
Since $w$ and $ut$ are injective, $d(ut)=d(t)$ and $d(t)=d(w)$, it follows that 
$w\in \genset{ut, \sym}\leq \genset{u, \all}.$ 
To summarise, $f \in \genset{\all, wv}\leq \genset{\all, u,v}$, as required.
\qed\vspace{\baselineskip}

We also require a dual of Lemma \ref{tech}, where the values of $c$ and $d$ are interchanged. 


\begin{lem} \label{tech2}
Let $u,v,f\in \T$. Then the following hold:
\begin{enumerate} [label = \rm{(\roman*})]
\item \label{d finite}
if $0<d(u), d(f)<\infty$, $c(u)=c(f)$ and $k(f) < \infty$, then   $f\in \genset{\all, u}$;
\item \label{d finite more}
if $c(u)=\infty$, $d(u)<\infty$, $c(f)=c(v)=0$, $d(v)>0$, and $0<d(f)<\infty$, then 
$f\in \genset{\all, u,v}$;
\item \label{c infinite} if $u$ is surjective, $d(f), k(f)<\infty$ and $c(u)=c(f)=\infty$, 
then $f\in \genset{\all, u}$;
\item \label{two fun dual} if $u$ is surjective, $c(u) > 0$, $c(v)=c(f)=\infty$, $k(f)< \infty$ and 
$d(v), d(f)<\infty$, then $f \in \genset{\all, u, v}$.
\end{enumerate}
\end{lem}
\proof 
{\bf (i).}   Let $g,h \in \T$ be such that $\Omega u$ is a transversal of
$g$, $d(g)=d(f)$, $\ker(h)=\ker(f)$, and $\Omega h$ is a transversal 
of $u$. Then $\ker(f)=\ker(hug)$ and 
$d(f)=d(hug)$. Hence $f\in \genset{hug, \sym}$ and so it suffices to show that $g,h\in \all$. 

Since $c(g)=d(u)$ and $d(g)=d(f)$,  it follows that $0<c(g),d(g)<\infty$ and so $g \in \all$.
Also $c(h)=c(f)$, $d(h)=c(u)$, and $k(h)=k(f)<\infty$. So, if $c(f)=c(u)=0$, then $h\in \sym$. If  $0<c(f)=c(u)<\infty$, then $0<c(h), d(h)<\infty$ and thus $h \in \all$. Finally, if $c(f)=c(u)=\infty$, then  $c(h)=d(h)=\infty$, $k(h)<\infty$, 
and so $h\in\all$.\vspace{\baselineskip}

\noindent {\bf (ii).} If $d(v)<\infty$, then $f \in \genset{\all, v}$ by Lemma \ref{tech2}\ref{d finite}. Suppose $d(v)=\infty$. 
Let $g \in \T$ be such that $\Omega v$ is a transversal of $g$, $\Omega g$ is contained in a transversal $\Sigma$ of $u$ with 
$0<|\Sigma\setminus \Omega g|<\infty$, and $k(g)<\infty$. 
Then $c(g)=d(g)=\infty$ and so $g \in \all$. Also $c(vgu)=0$ and 
$0<d(vgu)<\infty$. Therefore $f \in \genset{\all, vgu} \leq \genset{\all,u, v}$ by Lemma \ref{tech2}\ref{d finite}.
\vspace{\baselineskip}

\noindent {\bf (iii).} Let $g \in \T$ be such that $\ker(g)=\ker(f)$ and $\Omega g$ is contained in a
transversal $\Sigma$ of $u$ with $|\Sigma\setminus \Omega g| = d(f)$. Then $c(g)=d(g)=\infty$, $k(g)<\infty$ and 
thus $g \in \all$.  Also $\ker(f)=\ker(gu)$ and $d(f)=d(gu)$, and so $f\in \genset{gu, \sym}\leq \genset{u, \all}$.\vspace{\baselineskip}

\noindent {\bf (iv).} If $d(v)=0$ or $c(u)=\infty$, then the result follows from \ref{c infinite}. 
Suppose that $0<c(u),d(v)<\infty$ and let $g \in \T$ be such that $\Omega v$ is a transversal 
of $g$ and $\Omega g$ is a transversal of $u$. Then $c(g)=d(v)$ and $d(g)=c(u)$ which implies $0<c(g), d(g)<\infty$
and $g \in \all$. Therefore $f \in \genset{\all, vgu}$ by \ref{c infinite} since $c(vgu)=\infty$, 
$k(vgu)<\infty$ and $d(vgu)=0$. \qed \vspace{\baselineskip}


The final technical lemma we require relates to generating mappings with infinite $k$ value, whereas the previous two lemmas 
are concerned with generating mappings with finite $k$ value.

\begin{lem} \label{tech k}
Let $u,v,t,f\in \T$. Then the following hold:
\begin{enumerate} [label = \rm{(\roman*})]
\item\label{kd infinite} if $k(u)=k(f)=d(f)=\infty$, then $f\in \genset{\all, u}$;
\item\label{k d>0} if $k(u)=k(f)=\infty$, $d(f),d(u)>0$, $c(v)=\infty$, and $0<d(v)<\infty$ then
$f \in \genset{\all, u, v}$;
\item\label{k infinite} if $k(u)=k(f)=\infty$, $c(v)=\infty$, and $d(v)=0$, then $f \in \genset{\all,u,v}$;
\item\label{k three} if $k(u)=k(f)=\infty$, $c(v)=\infty$, $d(v)<\infty$, $c(t)>0$, and $d(t)=0$, then 
$f \in \genset{\all,u,v,t}$.
\end{enumerate}
\end{lem}
\proof 
{\bf (i).} Let $\set{K_i}{i\in \N}$ be  the  kernel classes of $f$ where $K_0$ is infinite.
Also let $\set{L_i}{i\in \N}$ be  the  infinite kernel classes of $u$. 
Let $g \in \Omega^\Omega$ be such that $K_0g=\{\alpha\}$ where $\alpha \in L_0$, $g|_{K_i}$ 
is injective and $K_ig\subseteq L_{2i}$. Then $c(g)=d(g)=\infty$ and $k(g)=1$ which implies that  
$g\in \all$. Also $\ker(gu)=\ker(f)$ and 
$d(gu)=d(f)=\infty$, and so $f\in \genset{gu, \sym}\leq \genset{u, \all}$. 
\vspace{\baselineskip}

\noindent{\bf (ii).} If $d(f) = \infty$, then the result follows from \ref{kd infinite}. 
Suppose $0<d(f)<\infty$ and let $g\in\Omega^\Omega$ be such that 
$\ker(f)=\ker(g)$ and $\Omega g$ is a transversal of $v$. 
Then $k(g)=k(f)=\infty$ and $d(g)=c(v)=\infty$ and so by \ref{kd infinite}, 
$g \in \genset{\all, u}$. 
If $h \in \T$ is such that $\Omega v$ is a transversal of $h$ and $d(h)=d(f)$, then 
$0<c(h), d(h)<\infty$ and $h\in \all$.
 Hence  $\ker(f)=\ker(gvh)$ and $d(f)=d(gvh)$, which implies that 
$f\in \genset{g,h, v, \sym}\leq \genset{u,v, \all}$.
\vspace{\baselineskip}

\noindent{\bf (iii).} Let $g \in \T$ be such that $\ker(f)=\ker(g)$ and $\Omega g$ is contained in a 
transversal 
$\Sigma$ of $u$ with $|\Sigma \setminus \Omega g| = d(f)$. Then $k(g)=d(g)=\infty$ and
$g \in \genset{\all, u}$  by \ref{kd infinite}. Since $\ker(f)=\ker(gu)$ and $d(f)=d(gu)$, 
$f \in \genset{gu, \sym}\leq \genset{\all, u}$, as required.
\vspace{\baselineskip}

\noindent{\bf (iv).} If $d(v)=0$ or $c(t)=\infty$, then $f\in \genset{u,v, \all}$ or $f\in \genset{u,t,\all}$, respectively, from \ref{k infinite}. Suppose
$0<c(t),d(v)<\infty$ and let $g \in \T$ be such that $\Omega v$ is a transversal 
of $g$ and $\Omega g$ is a transversal of $t$. Then $0<c(g), d(g)<c(t)$ and so $g\in \all$.  Also 
$c(vgt)=\infty$ and $d(vgt)=0$, and hence
$f \in \genset{\all, u, vgt}\leq \genset{u,v,t,\all}$ by \ref{k infinite}.
\qed


\section{The Proof of the Main Theorem} \label{sec:proof}

In this section we prove the Main Theorem.  We start by noting that since the intersection of two 
subsemigroups of a 
semigroup is a subsemigroup, it follows that the intersection of any of $S_1$, $S_2$, $S_3$,  
$S_4$ and  $S_5$ is 
subsemigroup of $\Omega^\Omega$. In the following proposition we prove that $U$ and $V$ 
are semigroups.


\begin{prop}\label{U_are_a_semigroup}
$U=\{ f \in \T : d(f) = \infty \text{ or } 0<c(f)<\infty \}\cup \sym$ and $V=  \{f \in \T : c(f)=\infty \text{ or } (0<d(f)<\infty) \} \cup \sym$  are 
semigroups.
\end{prop}
\proof 
Let $f,g\in U$. If $f\in \sym$, then $d(fg)=d(g)$ and $c(fg)=c(g)$ and so $fg\in U.$ Similarly, if $g\in \sym$, then $fg\in U$. 
If $d(g)=\infty$, then $d(fg)=\infty$ and $fg\in U$. If $d(f)=\infty$ and $0<c(g)<\infty$, then by Lemma \ref{inequalities}(iii), 
$d(fg)=\infty$ and so $fg\in U$.  If  $0<c(f), c(g)<\infty$, then, by Lemma \ref{inequalities}(ii), $0<c(fg)<\infty$ and so 
$fg\in U$. 

The proof that $V$ is a semigroup follows by a similar argument using Lemma \ref{inequalities}(i) and (iv).
\qed\vspace{\baselineskip}
 
 
We prove the Main Theorem by finding the maximal subsemigroups in $S_1$, $S_2$, $S_3$,  $S_4$ and  $S_5$ containing $\all$, then the 
maximal subsemigroups in each of  those 
semigroups which contain $\all$ and so on. We give the descriptions of the maximal subsemigroups found at each stage in a separate 
statement, each of which can be proved using the following strategy. 

Let $S$ be a semigroup, let $T$ be a subsemigroup of $S$, and let 
$\mathcal{M}=\set{M_i}{i\in I}$ be a collection of subsemigroups of $S$ containing $T$ for some set $I$ and such that 
 $M_i\not\leq M_j$ for all $i,j\in I$ with $i\not=j$. 
Suppose that $\mathcal{M}$ has the following property: 
\begin{quote}
if $U$ is a subsemigroup of $S$ containing $T$ and $U$ is not contained in any  $M_i$ for any $i\in I$, then $U=S$.
\end{quote}
Then it is routine to verify that $\mathcal{M}$ consists of the maximal subsemigroups of $S$ containing $T$.

There are essentially 36 cases in the proof of the Main Theorem. However, there are 14 pairs of cases where the proof of one case can be
obtained from the proof of the other by interchanging the values of $c$ and $d$, interchanging any part of Lemma \ref{tech} by the same 
part of Lemma \ref{tech2} (or vice versa), and by interchanging part (ii), (iii), or (iv) of Lemma \ref{tech k} with part (i) of the same lemma. 
Therefore we will not present duplicate proofs, but will only give a proof of one of the cases in each of these pairs.


\begin{prop} \mbox{}
\begin{enumerate}[label = \emph{(\roman*)}] 
  \item  the maximal subsemigroups of  $S_1$  containing $\all$ are: $S_{1,2}$, $S_{1,3}$, $S_{1,4}$ and $S_{1,5}$;
  \item the maximal subsemigroups of $S_2$ containing $\all$ are  $S_{1,2} $, $S_{1,3}$, $S_{1,4}$ and $S_{1,5}$;
   \item the maximal subsemigroups of $S_3$ containing $\all$ are  $S_{1,3}$, $U$, $S_{3,4}$ and $S_{3,5}$;
    \item  the maximal subsemigroups of $S_4$ containing $\all$  are: $V $, $S_{2,4}$, $S_{3,4}$ and $S_{4,5}$;
  \item the maximal subsemigroups of $S_5$ containing $\all$ are: $S_{1,5}$, $S_{2,5}$, $S_{3,5}$ and $S_{4,5}$. 
\end{enumerate}
\end{prop}
 \proof
 {\bf (i).}  
  It suffices to show that if $M$ is any subsemigroup of $S_1$ containing $\all$ but not 
  contained in any of the semigroups $S_{1,2}$, $S_{1,3}$, $S_{1,4}$ and $S_{1,5}$, then $M=S_1$. 
    Let $u_1 \in M \setminus S_{1,2}$, $u_2 \in M \setminus S_{1,3}$, $u_3 \in M \setminus 
    S_{1,4}$ and $u_4 \in M \setminus S_{1,5}$. Then the following hold:
\begin{align*}
 c(u_1)=0 \and d(u_1) >0, &\qquad c(u_2)= \infty \and  0 < d(u_2) < \infty\\
 c(u_3) < \infty \and d(u_3) = \infty,& \qquad k(u_4) = \infty \and d(u_4)>0.
\end{align*}
  Let $f\in S_1\setminus \all$ be arbitrary. Then one of the following holds: 
  \begin{enumerate}

    \item $f$ is injective and $0<d(f)<\infty$ in which case by Lemma \ref{tech2}\ref{d finite more}, $f\in \genset{\all, u_1, u_2}$;
    \item $0<d(f)<\infty$, $k(f)<\infty$ and $c(f)=\infty$ in which case by Lemma \ref{tech2}\ref{d finite}, $f\in \genset{\all, u_2}$;
    \item $c(f)<\infty$ and $d(f)=\infty$  in which case by Lemma \ref{tech}\ref{two fun}, $f\in \genset{\all, u_1, u_3}$;
    \item $k(f)=\infty$ and $d(f)>0$  in which case by Lemma \ref{tech k}\ref{k d>0}, $f\in \genset{\all, u_2, u_4}$. 
    \end{enumerate}
Therefore
$$S_1\leq \genset{\all, u_1, u_2, u_3, u_4}\leq M\leq S_1,$$
giving equality throughout.  \vspace{\baselineskip}
  
\noindent {\bf (ii).} 
  Let  $M$ be any subsemigroup of $S_2$ containing $\all$ which is not contained in any of the
  given semigroups. 
  Let $u_1 \in M \setminus S_{1,2}$, $u_2 \in M \setminus S_{2,3}$, $u_3 \in M \setminus 
    S_{2,4}$ and $u_4 \in M \setminus S_{2,5}$. Then
\begin{align*}
 c(u_1)>0 \and d(u_1) =0, &\qquad c(u_2)= \infty \and  d(u_2) < \infty\\
 0<c(u_3) < \infty \and d(u_3) = \infty,& \qquad k(u_4) = \infty.
\end{align*}
  Let $f\in S_2\setminus \all$ be arbitrary. Then one of the following holds: 
  \begin{enumerate}
    \item $0<c(f)<\infty$ and $d(f)=0$ in which case by Lemma \ref{tech}\ref{c finite more}, $f\in \genset{\all, u_1,u_3}$;
    \item $0<c(f)<\infty$ and $d(f)=\infty$ in which case by Lemma \ref{tech}\ref{c finite}, $f\in \genset{\all, u_3}$;
    \item $d(f)<\infty$, $k(f)<\infty$ and $c(f)=\infty$  in which case by Lemma \ref{tech2}\ref{two fun dual}, $f\in \genset{\all, u_1, u_2}$;
    \item $k(f)=\infty$  in which case by Lemma \ref{tech k}\ref{k three}, $f\in \genset{\all, u_1, u_2, u_4}$. 
    \end{enumerate}
Therefore
$$S_2\leq \genset{\all, u_1, u_2, u_3, u_4}\leq M\leq S_2,$$
giving equality throughout.  \vspace{\baselineskip}
  
\noindent {\bf (iii). }  The proof of this case follows by an argument analogous to that used in the proof of part (iv) as discussed before the 
proposition. It is also necessary in case (4) of (iv) to replace the assumption that $k(f)=\infty$  by  $k(f)=\infty$ 
and $d(f)=\infty$ and to apply Lemma \ref{tech k}\ref{kd infinite}. \vspace{\baselineskip}

\noindent {\bf (iv).}  
  Let  $M$ be any subsemigroup of $S_4$ containing $\all$ which is not contained in any of the
 given semigroups.
  Let $u_1 \in M \setminus V$, $u_2 \in M \setminus S_{2,4}$, $u_3 \in M \setminus S_{3,4}$ and $u_4 \in M \setminus S_{4,5}$. Then
\begin{align*}
 0<c(u_1)< \infty\and d(u_1) =0, &\qquad c(u_2)= 0 \and  0 < d(u_2) < \infty\\
 c(u_3) = \infty\text{ and  }d(u_3) < \infty,& \qquad k(u_4) = \infty.
\end{align*}

  Let $f\in S_4\setminus \all$ be arbitrary. Then one of the following holds: 
\begin{enumerate}
  \item $d(f)<\infty$, $k(f)<\infty$ and $c(f)=\infty$  in which case by Lemma \ref{tech2}\ref{two fun dual}, $f\in \genset{\all, u_1, u_3}$;
  \item $c(f)=0$ and $0<d(f)<\infty$ in which case by Lemma \ref{tech2}\ref{d finite}, $f\in \genset{\all, u_2}$;
  \item $0<c(f)<\infty$ and $d(f)=0$ in which case by Lemma \ref{tech}\ref{c finite}, $f\in \genset{\all, u_1}$;
    \item $k(f)=\infty$ in which case by Lemma \ref{tech k}\ref{k three}, $f\in \genset{\all, u_1, u_3, u_4}$.
  \end{enumerate}
Therefore
$$S_4\leq \genset{\all, u_1, u_2, u_3, u_4}\leq M\leq S_4,$$
giving equality throughout. \vspace{\baselineskip}

\noindent {\bf (v).} 
  Let  $M$ be any subsemigroup of $S_5$ containing $\all$ which is not contained in any of the
  given semigroups.
  Let $u_1 \in M \setminus S_{1,5}$, $u_2 \in M \setminus S_{2,5}$, $u_3 \in M \setminus S_{3,5}$ and $u_4 \in M \setminus S_{4,5}$. Then
\begin{align*}
 c(u_1)>0, k(u_1)<\infty \and d(u_1) =0, &\qquad c(u_2)= 0 \and  d(u_2)>0\\
 c(u_3) = \infty, k(u_3)<\infty \and d(u_3) < \infty,& \qquad c(u_4)<\infty \and d(u_4) = \infty.
\end{align*}

Let $f\in S_5\setminus \all$ be arbitrary. Then one of the following holds: 
\begin{enumerate}
  \item $d(f)<\infty$, $k(f)<\infty$ and $c(f)=\infty$  in which case by Lemma \ref{tech2}\ref{two fun dual}, $f\in \genset{\all, u_1, u_3}$;
  \item $c(f)<\infty$ and $d(f)=\infty$ in which case by Lemma \ref{tech}\ref{two fun}, $f \in \genset{\all, u_2, u_4}$;
  \item $0<c(f)<\infty$ and $d(f)=0$ in which case by Lemma \ref{tech}\ref{c finite more}, $f \in \genset{\all, u_1 ,u_4}$;
  \item $c(f)=0$ and $0<d(f)<\infty$ in which case by Lemma \ref{tech2}\ref{d finite more}, $f \in \genset{\all, u_2 ,u_3}$.
  \end{enumerate}
Therefore
$$S_5\leq \genset{\all, u_1, u_2, u_3, u_4}\leq M\leq S_5,$$
giving equality throughout.  \qed


\begin{prop} \mbox{}
\begin{enumerate}[label = \emph{(\roman*)}] 
  \item  the maximal subsemigroups of  $V$  containing $\all$ are: $S_{1,4}$, $V \cap S_2$ and $V \cap S_5$;
  \item the maximal subsemigroups of $U$ containing $\all$ are: $U\cap S_1$, $S_{2,3}$ and $U \cap S_5$.
\end{enumerate}
\end{prop}
\proof {\bf (i).}
  Let $u_1 \in V \setminus S_{1,4}$, $u_2 \in V \setminus (V \cap S_2)$ and $u_3 \in V \setminus (V \cap S_5)$. Then
  $$c(u_1)= \infty \and d(u_1) =0, \quad c(u_2)= 0 \and  0<d(u_2) <\infty, \quad k(u_3) = \infty.$$
  Let $f\in V\setminus \all$ be arbitrary. Then one of the following holds: 
\begin{enumerate}
  \item $0<d(f)<\infty$ and $c(f)=0$ in which case by Lemma \ref{tech2}\ref{d finite}, $f\in \genset{\all, u_2}$;
  \item $d(f)<\infty$, $k(f)<\infty$ and $c(f)=\infty$ in which case by Lemma \ref{tech2}\ref{c infinite}, $f\in \genset{\all, u_1}$;
  \item $k(f)=\infty$ in which case by Lemma \ref{tech k}\ref{k infinite}, $f \in \genset{\all, u_1, u_3}$. 
\end{enumerate}
Hence if $M$ is any subsemigroup of $V$ containing $\all$ which is not contained in any of the semigroups in the statement of the propoposition, then $M=V$. \vspace{\baselineskip}

\noindent {\bf (ii).} The proof is analogous to (i).
\qed


\begin{prop} \mbox{}
\begin{enumerate}[label = \emph{(\roman*)}] 
  \item the maximal subsemigroups of  $S_{1,2}$ containing $\all$ are: $S_{1,2,3}$, $S_{1,2,4}$ and $S_{1,2,5}$;
  \item the maximal subsemigroups of  $S_{1,3}$  containing $\all$ are: $U\cap S_1$, $S_{1,3,4}$ and $S_{1,3,5}$;
  \item the maximal subsemigroups of  $S_{1,4}$ containing $\all$ are: $S_{1,2,4}$, $S_{1,3,4}$ and $S_{1,4,5}$;
  \item the maximal subsemigroups of  $S_{1,5}$ containing $\all$ are: $S_{1,2,5}$, $S_{1,3,5}$ and $S_{1,4,5}$;
  \item the maximal subsemigroups of  $S_{2,3}$ containing $\all$ are: $S_{1,2,3} $, $S_{2,3,4}$ and $S_{2,3,5}$;
  \item the maximal subsemigroups of  $S_{2,4}$ containing $\all$ are: $V\cap S_2$, $S_{2,3,4}$ and $S_{2,4,5}$;
  \item the maximal subsemigroups of  $S_{2,5}$ containing $\all$ are: $S_{1,2,5}$, $S_{2,3,5}$ and $S_{2,4,5}$;
  \item the maximal subsemigroups of  $S_{3,4}$ containing $\all$ are: $S_{1,3,4}$, $S_{2,3,4}$ and $S_{3,4,5}$;
  \item the maximal subsemigroups of  $S_{3,5}$ containing $\all$ are: $S_{1,3,5}$, $U \cap S_5$ and $S_{3,4,5}$;
  \item the maximal subsemigroups of  $S_{4,5}$ containing $\all$ are: $V \cap S_5$, $S_{2,4,5}$ and $S_{3,4,5}$.
\end{enumerate}
\end{prop}
\proof
{\bf (i).} 
Let  $M$ be any subsemigroup of $S_{1,2}$ containing $\all$ which is not contained in any of 
the given semigroups. Then there exist $u_1, u_2, u_3\in M$ such 
that:
 $$c(u_1)=\infty \and 0< d(u_1) < \infty,\quad 0<c(u_2) < \infty \and d(u_2) = \infty, \quad k(u_3) = \infty \and  d(u_3) > 0.$$
 
 Let $f\in S_{1,2} \setminus \all$ be arbitrary. Then one of the following holds: 
\begin{enumerate}
\item  $d(f)=\infty$ and $0<c(f)<\infty$ in which case  Lemma \ref{tech}\ref{c finite}, $f\in \genset{\all, u_2}$;
\item $c(f)=\infty$, $k(f)<\infty$ and $0<d(f)<\infty$ in which case  Lemma \ref{tech2}\ref{d finite}, $f\in \genset{\all, u_1}$;
\item $k(f)=\infty$ and $d(f)>0$ in which case  Lemma \ref{tech k}\ref{k d>0}, $f\in \genset{\all, u_1, u_3}$.
\end{enumerate}
 Thus $M=S_{1,2}$.\vspace{\baselineskip}

\noindent{\bf (ii).} The proof is analogous to (vi). \vspace{\baselineskip}
 
\noindent{\bf (iii).} 
  Let  $M$ be any subsemigroup of $S_{1,4}$ containing $\all$ which is not contained in any of the
  given semigroups. Then there exist $u_1, u_2, u_3\in M$ such 
  that:
 $$c(u_1)=0 \and  0<d(u_1)< \infty, \quad c(u_2)=\infty  \and 0< d(u_2)< \infty, \quad 
 k(u_3) = \infty  \and d(u_3) >0.$$
 
 Let $f\in S_{1,4} \setminus \all$ be arbitrary. Then one of the following holds:
\begin{enumerate}
\item  $0<d(f)<\infty \and c(f)=0$ and so $f\in \genset{\all, u_1}$ by Lemma \ref{tech2}\ref{d finite};
\item  $0<d(f)<\infty$, $k(f)<\infty$ and $c(f)=\infty$ and so $f\in \genset{\all, u_2}$ by Lemma \ref{tech2}\ref{d finite};
\item $k(f)=\infty$ and $d(f)>0$ and so $f\in \genset{\all, u_2, u_3}$ by Lemma \ref{tech k}\ref{k d>0}.
\end{enumerate}
 Thus $M=S_{1,4}$, as required.\vspace{\baselineskip}
 
\noindent{\bf (iv).} 
  Let  $M$ be any subsemigroup of $S_{1,5}$ containing $\all$ which is not contained in any of the
given semigroups. Then there exist $u_1, u_2, u_3\in M$ such
that:
 $$c(u_1)=0 \and d(u_1)>0$$ 
 $$c(u_2)=\infty,\quad k(u_2)<\infty  \and 0< d(u_2)< \infty, \quad 
 d(u_3) = \infty  \and c(u_3) <\infty.$$
 
  Let $f\in S_{1,5} \setminus \all$ be arbitrary. Then one of the following holds:
\begin{enumerate}
\item  $0<d(f)<\infty \and c(f)=0$ and so $f\in \genset{\all, u_1, u_2}$ by Lemma \ref{tech2}\ref{d finite more};
\item  $0<d(f)<\infty$, $k(f)<\infty$ and $c(f)=\infty$ and so $f\in \genset{\all, u_2}$ by Lemma \ref{tech2}\ref{d finite};
\item $d(f)=\infty$ and $c(f)<\infty$ and so $f\in \genset{\all, u_1, u_3}$ by Lemma \ref{tech}\ref{two fun}.
\end{enumerate}
 Thus $M=S_{1,5}$, as required.\vspace{\baselineskip}
  
\noindent{\bf (v).} The proof is analogous to (iii). \vspace{\baselineskip}
 
 \noindent{\bf (vi).}  
 Let  $M$ be any subsemigroup of $S_{2,4}$ containing $\all$ which is not contained in any of 
 the given semigroups. Then there exist $u_1, u_2, u_3\in M$ such
  that:
 $$d(u_1)=0 \and 0< c(u_1) < \infty,\quad d(u_2) < \infty \and c(u_2) = \infty, \quad k(u_3) = \infty .$$
 
 Let $f\in S_{2,4} \setminus \all$ be arbitrary. Then one of the following holds: 
\begin{enumerate}
\item  $d(f)=0$ and $0<c(f)<\infty$ in which case  Lemma \ref{tech}\ref{c finite}, $f\in \genset{\all, u_1}$;
\item $c(f)=\infty$, $k(f)<\infty$ and $d(f)<\infty$  in which case  Lemma \ref{tech2}\ref{two fun dual}, $f\in \genset{\all, u_1, u_2}$;
\item $k(f)=\infty$  in which case  Lemma \ref{tech k}\ref{k three}, $f\in \genset{\all, u_1, u_2, u_3}$.
\end{enumerate}
 Thus $M=S_{2,4}$.\vspace{\baselineskip} 
 
 \noindent{\bf (vii).}  The proof is analogous to (iv). \vspace{\baselineskip}
 
 \noindent {\bf (viii).}
Let  $M$ be any subsemigroup of $S_{3,4}$ containing $\all$ which is not contained in any of the
  given semigroups. Then there exist $u_1, u_2, u_3\in M$ such that:
$$0 < c(u_1) < \infty  \and  d(u_1)=0,\quad
    c(u_2)= 0  \and 0 < d(u_2) < \infty, \quad
    k(u_3) = \infty  \and  d(u_3) = \infty.$$
    
 Let $f\in S_{3,4} \setminus \all$ be arbitrary. Then one of the following holds:
\begin{enumerate}
\item $0<c(f)<\infty \and d(f)=0$ and so $f\in \genset{\all, u_1}$ by Lemma \ref{tech}\ref{c finite};
\item $c(f)=0 \and 0<d(f)<\infty$ and so $f\in \genset{\all, u_2}$ by Lemma  \ref{tech2}\ref{d finite};
\item $k(f)=d(f)=\infty$ and so $f\in  \genset{\all, u_3}$ by Lemma \ref{tech k}\ref{kd infinite}.
\end{enumerate}
Therefore $M=S_{3,4}$.\vspace{\baselineskip}
 
\noindent   {\bf (ix).} 
Let $M$ be any subsemigroup of $S_{3,5}$ containing $\all$ which is not contained in any of the
  given semigroups. Then there exist $u_1, u_2, u_3\in M$ such that:
$$0 < c(u_1) < \infty \and d(u_1) = 0, \quad c(u_2) < \infty\and  d(u_2) = \infty, \quad
   c(u_3) = 0 \and 0 < d(u_3) < \infty.$$

 Let $f\in S_{3,5}\setminus \all$ be arbitrary. Then one of the following holds:
\begin{enumerate}
\item $0<c(f)<\infty \and d(f)=0$ and so $f\in \genset{\all, u_1}$ by Lemma \ref{tech}\ref{c finite};
\item $c(f)<\infty \and d(f)=\infty$ and so $f\in \genset{\all, u_2, u_3}$ by Lemma \ref{tech}\ref{two fun};
\item $c(f)=0 \and 0<d(f)<\infty$ and so $f\in \genset{\all,  u_3}$ by Lemma \ref{tech2}\ref{d finite}.
\end{enumerate}
Whence $M=S_{3,5}$.\vspace{\baselineskip}

\noindent{\bf (x).} The proof is analogous to (ix). \qed


\begin{prop} \mbox{}
\begin{enumerate}[label = \emph{(\roman*)}] 
  \item the maximal subsemigroups  of $V\cap S_2$  containing $\all$ are: $S_{1,2,4}$ and $ V\cap S_{2,5}$;
  \item the maximal subsemigroups  of $U\cap S_1$  containing $\all$ are: $S_{1,2,3}$ and $ U \cap S_{1,5}$;
  \item the maximal subsemigroups  of $V\cap S_5$  containing $\all$ are: $S_{1,4,5}$ and $ V \cap S_{2,5}$;
  \item  the maximal subsemigroups of $ U\cap S_5$ containing $\all$ are: $S_{2,3,5}$ and $U \cap S_{1,5}$.
\end{enumerate}
\end{prop} 
\proof
{\bf (i).} 
Let $u_1\in (V\cap S_2)\setminus S_{1,2,4}$ and let  $u_2\in (V\cap S_2)\setminus (V\cap S_{2,5})$. Then:
 $$d(u_1)=0 \and c(u_1)=\infty, \quad k(u_2)=\infty.$$
  
  Let $f\in V\cap S_2 \setminus \all$ be arbitrary. Then one of the following holds:
\begin{enumerate}
\item $d(f)<\infty$, $k(f)<\infty$ and $c(f)=\infty$ and so $f\in \genset{\all, u_1}$ by Lemma \ref{tech2}\ref{c infinite};
\item  $k(f)=\infty$ and so  $f\in \genset{\all, u_1, u_2}$ by Lemma \ref{tech k}\ref{k infinite}.
\end{enumerate}
So, if $M$ is any subsemigroup of $V\cap S_2$ containing $\all$ which is not contained in any of the
  given semigroups, then $M=V\cap S_2.$\vspace{\baselineskip}

\noindent {\bf (ii).} The proof is analogous to (i).\vspace{\baselineskip}
  
\noindent   {\bf (iii).}  
  Let $u_1 \in (V \cap S_5)\setminus S_{1,4,5}$ and $u_2 \in (V \cap S_5) \setminus 
  (V\cap S_{2,5})$. 
  Then $$c(u_1)=\infty,\ k(u_1)<\infty \and  d(u_1)=0, \quad c(u_2)=0 \and 0<d(u_2)<\infty.$$
  
If $f\in (V\cap S_5)\setminus \all$ is arbitrary, then one of the following holds:
\begin{enumerate}
\item $c(f)=\infty$, $k(f)<\infty$ and $d(f)<\infty$ and so  Lemma \ref{tech2}\ref{c infinite} implies that $f\in \genset{\all, u_1}$;
\item $c(f)=0$, $0< d(f)<\infty$, and $k(f)<\infty$ and so Lemma \ref{tech2}\ref{d finite} implies that $f\in \genset{\all, u_2}$.
\end{enumerate}
So, if $M$ is any subsemigroup of $V\cap S_5$ containing $\all$ which is not contained in any of the
  given semigroups, then $M=V\cap S_5.$\vspace{\baselineskip}
  
\noindent {\bf (iv).} The proof is analogous to (iii). \qed


\begin{prop}\label{not_sure}\mbox{}
\begin{enumerate}[label = \emph{(\roman*)}] 
  \item the maximal subsemigroups  of $S_{1,2,3}$  containing $\all$ are: $S_{1,2,3,4}$ and $S_{1,2,3,5}$;
  \item the maximal subsemigroups  of $S_{1,2,4}$  containing $\all$ are: $S_{1,2,3,4}$ and $S_{1,2,4,5}$;
  \item the maximal subsemigroups  of $S_{1,2,5}$  containing $\all$ are: $S_{1,2,3,5}$ and $S_{1,2,4,5}$;
  \item the maximal subsemigroups  of $S_{1,3,4}$  containing $\all$ are: $S_{1,2,3,4}$ and $S_{1,3,4,5}$;
  \item the maximal subsemigroups  of $S_{1,3,5}$  containing $\all$ are: $S_{1,3,4,5}$ and $U\cap S_{1,5}$;
  \item the maximal subsemigroups  of $S_{1,4,5}$  containing $\all$ are: $S_{1,2,4,5}$ and $S_{1,3,4,5}$;
  \item the maximal subsemigroups  of $S_{2,3,4}$  containing $\all$ are: $S_{1,2,3,4}$ and $S_{2,3,4,5}$;
  \item the maximal subsemigroups  of $S_{2,3,5}$  containing $\all$ are: $S_{1,2,3,5}$ and $S_{2,3,4,5}$;
  \item the maximal subsemigroups  of $S_{2,4,5}$  containing $\all$ are: $S_{2,3,4,5}$ and $V \cap
  S_{2,5}$;
  \item the maximal subsemigroups  of $S_{3,4,5}$  containing $\all$ are: $S_{1,3,4,5}$ and $S_{2,3,4,5}$.
\end{enumerate}
\end{prop}
\proof
  {\bf (i).}  
  Let $u_1 \in S_{1,2,3}  \setminus S_{1,2,3,4}$ and let $u_2 \in S_{1,2,3}  \setminus (U\cap S_{1,5})$. Then 
$$ 0<c(u_1)<\infty \and  d(u_1)=\infty, \quad k(u_2)=d(u_2)=\infty.$$

 If $f\in S_{1,2,3}\setminus \all$ is arbitrary, then one of the following holds:
\begin{enumerate}
\item $0<c(f)<\infty$ and $d(f)=\infty$ and so  $f\in \langle u_1, \all \rangle$ by Lemma \ref{tech}\ref{c finite};
\item  $k(f)=d(f)=\infty$ and so $f\in \langle u_2, \all \rangle$ by Lemma \ref{tech k}\ref{kd infinite}. 
  \end{enumerate}
  So, if $M$ is any subsemigroup of $S_{1,2,3}$ containing $\all$ which is not contained in any of the
 given semigroups, then $M=S_{1,2,3}$.\vspace{\baselineskip}


  \noindent{\bf (ii).} 
  Let $u_1 \in S_{1,2,4}\setminus S_{1,2,3,4}$ and $u_2 \in S_{1,2,4}  \setminus S_{1,2,4,5}$. Then 
  $$c(u_1)=\infty,\ k(u_1)<\infty \and 0<d(u_1)<\infty,\quad k(u_2)=\infty \and d(u_2)>0.$$
   
 If $f\in S_{1,2,4}\setminus \all$ is arbitrary, then one of the following holds:
\begin{enumerate}
\item  $c(f)=\infty$, $k(f)<\infty$ and $0<d(f)<\infty$ and so $f \in \langle u_1, \all \rangle$ by Lemma \ref{tech2}\ref{d finite};
\item $k(f)=\infty$ and $d(f)>0$ and so $f \in \langle u_1, u_2, \all \rangle$ by Lemma \ref{tech k}\ref{k d>0}.
\end{enumerate}
So, if $M$ is any subsemigroup of $S_{1,2,4}$ containing $\all$ which is not contained in any of the
  given semigroups, then $M=S_{1,2,4}$.\vspace{\baselineskip}


  \noindent{\bf (iii).}
  Let $u_1 \in S_{1,2,5}\setminus S_{1,2,3,5}$ and $u_2 \in S_{1,2,5}  \setminus S_{1,2,4,5}$. Then 
  $$c(u_1)=\infty,\ k(u_1)<\infty \and 0<d(u_1)<\infty,\quad 0<c(u_2)<\infty \and d(u_2)=\infty.$$
   
 If $f\in S_{1,2,5}\setminus \all$ is arbitrary, then one of the following holds:
\begin{enumerate}
\item  $0<c(f)<\infty$ and $d(f)=\infty$ and so $f \in \langle u_2, \all \rangle$ by Lemma \ref{tech}\ref{c finite};
\item $c(f)=\infty$, $k(f)<\infty$ and $0<d(f)<\infty$ and so $f \in \langle u_1, \all \rangle$ by Lemma \ref{tech2}\ref{d finite}.
\end{enumerate}
So, if $M$ is any subsemigroup of $S_{1,2,5}$ containing $\all$ which is not contained in any of the
 given semigroups, then $M=S_{1,2,5}$.\vspace{\baselineskip}


  \noindent{\bf (iv).}  
  Let $u_1 \in S_{1,3,4}\setminus S_{1,2,3,4}$ and $u_2 \in S_{1,3,4}  \setminus S_{1,3,4,5}$. Then 
  $$c(u_1)=0\and 0<d(u_1)<\infty,\quad k(u_2)=d(u_2)=\infty.$$ 
 If $f\in S_{1,3,4}\setminus \all$ is arbitrary, then one of the following holds:
\begin{enumerate}
\item  $c(f)=0 \and 0<d(f)<\infty$ and so $f\in \langle u_1, \all \rangle$ by Lemma \ref{tech2}\ref{d finite};
\item $k(f)=d(f)=\infty$ and so $f\in \langle u_2, \all \rangle$ by Lemma \ref{tech k}\ref{kd infinite}.
\end{enumerate}
So, if $M$ is any subsemigroup of $S_{1,3,4}$ containing $\all$ which is not contained in any of the
  given semigroups, then $M=S_{1,3,4}$.\vspace{\baselineskip}

  
\noindent  {\bf (v).} 
  Let $u_1 \in S_{1,3,5}\setminus (U \cap S_{1,5})$ and $u_2 \in S_{1,3,5}\setminus  S_{1,3,4,5}$. Then 
$$c(u_1)=0 \and 0<d(u_1)<\infty, \quad c(u_2)<\infty \and  d(u_2)=\infty.$$ 
 If $f\in S_{1,3,5}\setminus \all$ is arbitrary, then one of the following holds:
\begin{enumerate}
\item $c(f)=0$ and $0<d(f)<\infty$ and so $f\in  \langle u_1, \all \rangle$ by Lemma \ref{tech2}\ref{d finite};
\item $c(f)< \infty$ and $d(f)=\infty$ and so $f\in \langle u_1, u_2, \all \rangle$ by Lemma \ref{tech}\ref{two fun}. 
\end{enumerate}
So, if $M$ is any subsemigroup of $S_{1,3,5}$ containing $\all$ which is not contained in any of the
 given semigroups, then $M=S_{1,3,5}$.\vspace{\baselineskip}

  
\noindent  {\bf (vi).} The proof is analogous to (viii).  \vspace{\baselineskip}
  
  
\noindent  {\bf (vii).}   The proof is analogous to (iv). \vspace{\baselineskip}
  

 \noindent {\bf (viii).}  
  Let $u_1 \in S_{2,3,5}\setminus S_{1,2,3,5}$ and $u_2 \in S_{2,3,5}\setminus S_{2,3,4,5}$. Then 
$$0<c(u_1)<\infty \and  d(u_1)=0, \quad 0<c(u_2)<\infty  \and d(u_2)=\infty.$$

 If $f\in S_{2,3,5}\setminus \all$ is arbitrary, then one of the following holds:
\begin{enumerate}
\item $0<c(f)<\infty \and d(f)=0$ and so $f\in  \langle u_1, \all \rangle$ by Lemma \ref{tech}\ref{c finite};
\item $0<c(f)<\infty \and d(f)=\infty$ and so $f\in \langle u_2, \all \rangle$ by Lemma \ref{tech}\ref{c finite}.
\end{enumerate}
So, if $M$ is any subsemigroup of $S_{2,3,5}$ containing $\all$ which is not contained in any of the
given semigroups, then $M=S_{2,3,5}$.\vspace{\baselineskip}


\noindent  {\bf (ix).}  The proof is analogous to (v).\vspace{\baselineskip}
 
  
\noindent  {\bf (x).} 
Let $u_1 \in S_{3,4,5} \setminus S_{1,3,4,5}$ and let $u_2 \in S_{3,4,5} \setminus S_{2,3,4,5}$. Then 
 $$0<c(u_1)<\infty \and d(u_1)=0, \quad c(u_2)=0 \and 0<d(u_2)<\infty.$$
 If $f\in S_{3,4,5}\setminus \all$ is arbitrary, then one of the following holds:
 \begin{enumerate}
\item $0<c(f)<\infty \and d(f)=0$ and $f\in \langle u_1, \all \rangle$ by Lemma \ref{tech}\ref{c finite};
\item  $c(f)=0 \and 0<d(f)<\infty$ and so $f\in \langle u_2, \all \rangle$ by Lemma \ref{tech2}\ref{d finite}. 
\end{enumerate}
So, if $M$ is any subsemigroup of $S_{3,4,5}$ containing $\all$ which is not contained in any of the
  given semigroups, then $M=S_{3,4,5}$.
\qed


\begin{prop}\mbox{}
\begin{enumerate}[label = \emph{(\roman*)}] 
  \item the only maximal subsemigroup of $V\cap S_{2,5}$ containing $\all$ is $S_{1,2,4,5}$.
  \item the only maximal subsemigroup of $U\cap S_{1,5}$ containing $\all$ is $S_{1,2,3,5}$.
\end{enumerate}
\end{prop}
\proof 
{\bf (i)}
The semigroup $S_{1,2,4,5}$ is described in Proposition \ref{not_sure}(ii). 
Let $u \in (V\cap S_{2,5}) \setminus S_{1,2,4,5}$. Then $c(u)=\infty$ and $d(u)=0$ (also $k(u)<\infty$ holds) and so
by Lemma \ref{tech2}\ref{c infinite} 
$$V\cap S_{2,5}=\all\cup\set{f\in \Omega^\Omega}{d(f)<\infty,\ k(f)<\infty, \and c(f)=\infty}\subseteq  \langle u, \all \rangle\subseteq \genset{u, S_{1,2,4,5}}.$$
Therefore  if $M$ is any subsemigroup of $V\cap S_{2,5}$ containing $\all$ which is not contained in $S_{1,2,4,5}$, then $M=V\cap S_{2,5}$.\vspace{\baselineskip}

\noindent{\bf (ii)} The proof is analogous to (i).
\qed


\begin{prop}
$\all$ is maximal in $S_{1,2,3,4}$, $S_{1,2,3,5}$, $S_{1,2,4,5}$, $S_{1,3,4,5}$, and 
$S_{2,3,4,5}$.
\end{prop}
\proof
If  $u \in S_{1,2,3,4}\setminus \all$, then $k(u)=d(u)=\infty$. Hence by Lemma \ref{tech k}\ref{kd 
infinite} 
 $$\langle u, \all \rangle  \supseteq \{ f \in \T : k(f)=d(f)=\infty \} \cup \all = S_{1,2,3,4},$$
and so $\all$ is maximal in $S_{1,2,3,4}$.
 
 
 If $u\in S_{1,2,3,5}\setminus \all$, then  $0<c(u)<\infty$ and $d(u)=\infty$. 
 Therefore  by Lemma \ref{tech}\ref{c finite} 
$$\langle u, \all \rangle  \supseteq \{ f \in \T : 0<c(f)<\infty \and d(f)=\infty \} \cup \all= S_{1,2,3,5}$$
and so $\all$ is maximal in $S_{1,2,3,5}$.
The proof that $\all$ is maximal in $S_{1,2,4,5}$ is analogous to the proof that $\all$ is maximal in $S_{1,2,3,5}$ using Lemma \ref{tech2}\ref{d finite}.
 

    
If $u \in S_{1,3,4,5}  \setminus \all$, then  $c(u)=0 \and 0<d(u)<\infty$. Thus 
by Lemma \ref{tech2}\ref{d finite}
    \begin{align*}
      \langle u, \all \rangle  \supseteq \{ f \in \T : c(f)=0 \and 0<d(f)<\infty \} \cup \all
                               = S_{1,3,4,5}.
    \end{align*}
    In particular,  $\all$ is maximal in $S_{1,3,4,5}$.
    The proof that $\all$ is maximal in $S_{2,3,4,5}$ is analogous to the proof that $\all$ is maximal in $S_{1,3,4,5}$ using Lemma \ref{tech}\ref{c finite}.\qed
    
    
  

\end{document}